\documentclass[12pt]{article}%

\nonstopmode

\usepackage[numbers]{natbib}
\usepackage{amsmath}
\usepackage{amsbsy}
\usepackage{amscd}
\usepackage{amsopn}
\usepackage{amstext}
\usepackage{amsxtra}
\usepackage{amssymb}
\usepackage{amsthm,amsfonts}
\usepackage{graphicx}
\usepackage{epstopdf}

\usepackage{epsfig}

\newtheorem{theorem}{Theorem}[section]
\newtheorem{proposition}[theorem]{Proposition}
\newtheorem{lemma}[theorem]{Lemma}

\newtheorem{corollary}[theorem]{Corollary}

\theoremstyle{definition}
\newtheorem{definition}{Definition}

\newtheorem{remark}{Remark}

\newtheorem{remarks}{Remarks}

\newtheorem{conjecture}{Conjecture}

\def\beq{\begin{eqnarray}}
\def\eeq{\end{eqnarray}}
\def\beqq{\begin{eqnarray*}}
\def\eeqq{\end{eqnarray*}}
\def\beeq{\begin{eqnarray*}}
\def\eeeq{\end{eqnarray*}}
\def\be{\begin{equation}}
\def\ee{\end{equation}}

\def\alphabf{{\boldsymbol\alpha}}

\def\BF={{\mathbb{F}}}
\def\BG={{\mathbb{G}}}
\def\BH={{\mathbb{H}}}
\def\bOne={{\bf 1}}

\def\CE{{\cal E}}

\def\CG{{\cal G}}

\def\CI{{\cal I}}

\def\CN{{\cal N}}

\def\pa{\partial}

\def\var{{\mbox{Var}}}

\def\qed{\hfill$\sqcap\kern-8.0pt\hbox{$\sqcup$}$\\}

\def\BBC{\mathbb{C}}
\def\BBE{\mathbb{E}}
\def\BBP{\mathbb{P}}

\def\BBR{\mathbb{R}}

\def\BBZ{\mathbb{Z}}

\def\indeg{\mbox{deg}^-}
\def\outdeg{\mbox{deg}^+}
\def\One{{\bf 1}}

\def\Parrow{\stackrel{\mathrm{P}}{\longrightarrow}}
\def\Pequals{\stackrel{\mathrm{P}}{=}}

\begin{document}

\title{The Construction and Properties of Assortative Configuration Graphs }
\author{T. R.  Hurd$^1$\\
\emph{ $^1$ Department of Mathematics, McMaster University, Canada }\\
 }
%
\date{October 17, 2015}
\maketitle

%
%
\maketitle

\abstract{In the new field of financial systemic risk, the network of interbank counterparty relationships can be described as a directed random graph. In ``cascade models'' of systemic risk, this ``skeleton'' acts as the medium through which financial contagion is propagated. It has been observed in real networks that such counterparty relationships exhibit {\it negative assortativity}, meaning that a bank's counterparties are more likely to have unlike characteristics. This paper introduces and studies a general class of random graphs  called the assortative configuration model, parameterized by an arbitrary node-type distribution $P$ and edge-type distribution $Q$. The first main result is a law of large numbers that says the empirical edge-type distributions converge in probability to $Q$. The second  main result is a  formula for the large $N$ asymptotic probability distribution of general graphical objects called ``configurations''. This formula exhibits a key property called ``locally tree-like'' that in simpler models is known to imply strong results of percolation theory on the size of large connected clusters. Thus this paper provides the essential foundations needed to prove rigorous percolation bounds and cascade mappings in assortative networks.   }

 
\bigskip The ``skeleton'' of a financial network at a moment in time is the directed graph whose directed edges indicate which pairs of banks are deemed to have a significant counterparty relationship at this time. The arrow on each edge points from debtor to creditor. 
 It has been often observed in financial networks (and as it happens, also the world wide web) that they are highly {\it disassortative}\index{disassortative}, or as we prefer to say, {\it negatively assortative} (see for example \cite{Soramakietal07} and \cite{BechAtal10}). This refers to the property that any bank's counterparties (i.e. their graph neighbours) have a marked tendency to be banks of an opposite character. For example, it is observed that small banks tend to lend preferentially to large banks rather than other small banks. On the other hand, social networks are commonly observed to have positive assortativity: highly popular people are more likely to have highly popular friends. Structural characteristics such as degree distribution and assortativity are felt by some (see \cite{MayArin10}) to be highly relevant to properties of networks, notably their susceptibility to the propagation of contagion effects. However, the nature of such relationships is far from clear. The present paper introduces and studies a class of assortative directed random graphs that is both rich enough to describe real financial, engineered and social networks, and amenable to analytic treatment.  In this class, one can hope to understand better the relations between the local network topology and percolation properties. 
  
The main aim of this paper is to put a firm theoretical foundation under the class of configuration graphs with arbitrary node degree distribution $P$ and edge degree distribution $Q$. The class of configuration graphs with general $Q$ has not been well studied previously, and we will generalize some of the classic large $N$ asymptotic results known to be true for the nonassortative configuration graph construction introduced by \cite{Bollobas80} and others, and described in Section \ref{Configuration}. At the end of the paper, a new approximate Monte Carlo simulation algorithm for assortative configuration graphs is proposed.

\section{Definitions and Basic Results}
\label{sec:1} This section provides some standard graph theoretic definitions and develops an efficient notation for what will follow. Since this paper deals only with directed graphs rather than undirected graphs,  the term ``graph'' will have that meaning. Undirected graphs fit in easily as a subcategory of the directed case. 

\begin{definition}
 \begin{enumerate}
  \item For any $N\ge 1$, the collection of {\it directed graphs}\index{graph!directed graph} on $N$ nodes is denoted $\CG(N)$. The set of {\it nodes}\index{node} $\CN$ is numbered by integers, i.e. $\CN=\{1,\dots, N\}:=[N]$. Then  $g\in\CG(N)$, a graph on $N$ nodes, is a pair $(\CN,\CE)$ where the set of edges is a subset  $\CE\subset\CN\times\CN$ and each element    $\ell\in\CE$  is an ordered pair $\ell=(v,w)$ called an {\it edge}\index{edge} or {\it link}\index{link}. Links are labelled by integers $\ell\in\{1, \dots, E\}:=[E]$ where $E=|\CE|$. Normally,  ``self-edges'' with $v=w$ are excluded from $\CE$, that is, $\CE\subset\CN\times\CN\setminus{\rm diag}$.
  \item A given graph $g\in\CG(N)$ can be represented by its  $N\times N$ {\it adjacency matrix} \index{adjacency matrix} $M(g)$ with components 
  \[ M_{vw}(g)=\left\{\begin{array}{ll }
  1    &  \mbox{if $(v,w)\in g$ } \\
   0   &    \mbox{if $(v,w)\in\CN\times\CN\setminus g$ } 
\end{array}\right.\ .
\]
\item The {\it in-degree}\index{in-degree} $\indeg(v)$ and {\it out-degree}\index{out-degree} $\outdeg(v)$ of a node $v$ are 
\[\indeg(v)=\sum_w M_{wv}(g),\quad \outdeg(v)=\sum_w M_{vw}(g)\ .\]
  \item A node $v\in\CN$ has {\it node type}\index{node type} $(j,k)$ if its in-degree is $\indeg(v)=j$ and its out-degree is $\outdeg(v)=k$; the node set partitions into node types, $\CN=\cup_{jk}\CN_{jk}$. One writes $k_v=k,j_v=j$ for any $v\in\CN_{jk}$ and allow degrees to be any non-negative integer. 
  \item An edge $\ell=(v,w)\in\CE=\cup_{kj}\CE_{kj}$ is said to have {\it edge type}\index{edge type} $(k,j)$ with in-degree $j$ and out-degree $k$ if it is an out-edge of a node $v$ with out-degree $k_v=k$ and an in-edge of a node $w$ with in-degree $j_w=j$. One writes $\outdeg(\ell)=k_\ell=k$ and $\indeg(\ell)=j_\ell=j$ whenever $ \ell\in\CE_{kj}$. 
\item For completeness,  an \index{graph!undirected graph} undirected graph is defined to be any directed graph $g$ for which $M(g)$ is symmetric. 
\end{enumerate}
\end{definition}

\begin{figure}
\includegraphics[scale=0.4]{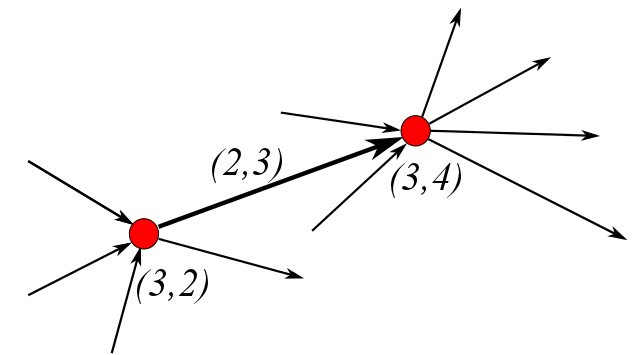}
 \caption{A type $(3,2)$ debtor bank that owes to a type $(3,4)$ creditor bank through a type $(2,3)$ link.}
\label{2nodes}
\end{figure}

The standard visualization of a graph $g$ on $N$ nodes is to plot nodes as ``dots'' with labels $v\in\CN$, and any edge $(v,w)$ as an arrow pointing ``downstream'' from node $v$ to node $w$. In the financial system application, such an arrow signifies that bank $v$ is a debtor of bank $w$ and the in-degree $\indeg(w)$ is the number of banks in debt to $w$, in other words the existence of the edge $(v,w)$ means ``$v$ owes $w$''.  Figure \ref{2nodes} illustrates the labelling of types of nodes and edges.

There are obviously constraints on the collections of node type $(j_v,k_v)_{v\in\CN}$ and edge type $(k_\ell,j_\ell)_{\ell\in\CE}$ if they derive from a graph. By computing the total number of edges $E=|\CE|$, the number of edges with $k_\ell=k$ and the number of edges with $j_\ell=j$, one finds three conditions: 
\beq
E:=|\CE|&=&\sum_v k_v=\sum_v j_v\nonumber\\
e^+_k:=|\CE\cap\{k_\ell=k\}| &=& \sum_\ell \One(k_\ell=k)=\sum_v k\One(k_v=k)\nonumber\\
\label{graphical}e^-_j:=|\CE\cap\{j_\ell=j\}| &=& \sum_\ell \One(j_\ell=j)=\sum_v j\One(j_v=j)\ .
\eeq

It is useful to define some further graph theoretic objects and notation in terms of the adjacency matrix $M(g)$:\begin{enumerate}
  \item The {\it in-neighbourhood}\index{in-neighbourhood} of a node $v$ is the set $\CN_v^-:=\{w\in\CN| M_{wv}(g)=1\}$ and the {\it out-neighbourhood}\index{out-neighbourhood}  of $v$ is the set $\CN_v^+:=\{w\in\CN| M_{vw}(g)=1\}$.
   \item One writes $\CE^+_v$  (or $\CE^-_v$) for the set of out-edges (respectively, in-edges) of a given node $v$ and $v^+_\ell$ (or $v^-_\ell$) for the node for which $\ell$ is an out-edge (respectively, in-edge). 
    \item Similarly, second-order neighbourhoods $\CN^{--}_v, \CN^{-+}_v, \CN^{+-}_v, \CN^{++}_v$ have the obvious definitions. Second and higher order neighbours can be determined directly from the powers of $M$ and $M^\top$. For example, $w\in \CN^{-+}_v$ whenever $(M^\top M)_{wv}\ge 1$. 
     \item One often writes  $j,j',j'',j_1$, etc. to refer to in-degrees and $k,k',k'',k_1$, etc. refer to out-degrees.
\end{enumerate}

Financial network models typically have a sparse adjacency matrix $M(g)$ when $N$ is large, meaning that the number of edges is a small fraction of the $N(N-1)$ potential edges. This reflects the fact that bank counterparty relationships are expensive to build and maintain, and thus $\CN^+_v$  and $\CN^-_v$ typically contain relatively few nodes even in a very large network. 

\subsection{Random Graphs}
Random graphs are simply probability distributions on the sets $\CG(N)$: 
\begin{definition}
 \begin{enumerate}
  \item A {random graph}\index{random graph} of size $N$ is a probability distribution $\BBP$ on the finite set $\CG(N)$. When the size $N$ is itself random, the probability distribution $\BBP$ is on the countable infinite set $\CG:=\cup_N \CG(N)$. Normally, it is assumed that $\BBP$ is invariant under permutations of the $N$ node labels. 
  \item For random graphs, define the {\it node-type distribution}\index{node-type distribution} to have probabilities $P_{jk}=\BBP[v\in\CN_{jk}]$ and the {\it edge-type distribution}\index{edge-type distribution} to have probabilities $Q_{kj}=\BBP[\ell\in\CE_{kj}]$. 
\end{enumerate}\end{definition}

$P$ and $Q$ can be viewed as bivariate distributions on the natural numbers, with marginals $P^+_k=\sum_j P_{jk}, P^-_j=\sum_k P_{jk}$ and $Q^+_k=\sum_j Q_{kj}, Q^-_j=\sum_k Q_{kj}$. Edge and node type distributions cannot be chosen independently however, but must be consistent with the fact that they derive from actual graphs which is true if one imposes that equations \eqref{graphical} hold in expectation, that is,  $P$ and $Q$ are ``consistent'':
\beq
   \nonumber z&:=&\sum_k kP^+_k=\sum_jjP^-_j\\
  \label{Consistent1} Q^+_k &=& kP^+_k/z, \quad Q^-_j = jP^-_j/z\quad \forall k,j\ .
  \eeq
  \label{PQassumptions}
  
A number of random graph construction algorithms have been proposed in the literature, motivated by the desire to create families of graphs that match the types and measures of network topology that have been observed in nature and society. The present paper focusses on so-called configuration graphs. The textbook ``Random Graphs and Complex Networks'' by van der Hofstad \cite{vdHofstad14} provides a complete and up-to-date review of the entire subject. 

In the analysis to follow, asymptotic results are typically expressed in terms of convergence of random variables in probability, defined as:
\begin{definition}
A sequence $\{X_n\}_{n\ge 1}$ of random variables is said to {\it converge in probability}\index{convergence in probability} to a random variable $X$, written $\lim_{n\to\infty} X_n\Pequals X$ or $X_n\Parrow X$, if for any $\epsilon>0$ 
\[ \BBP[|X_n-X|>\epsilon]\to 0\ .
\] \end{definition}
Recall further standard notation for asymptotics of sequences of real numbers\\ $\{x_n\}_{n\ge 1}, \{y_n\}_{n\ge 1}$ and random variables $\{X_n\}_{n\ge 1}$: \begin{enumerate}
  \item Landau's ``little oh'':  $x_n=o(1)$  means  $x_n\to 0$; $x_n=o(y_n)$ means $x_n/y_n=o(1)$;
  \item Landau's ``big oh'':  $x_n=O(y_n)$ means there is $N>0$ such that $x_n/y_n$ is bounded for $n\ge N$; 
  \item $x_n\sim y_n$ means $x_n/y_n\to 1$;
  \item $X_n\Pequals o(y_n)$ means $X_n/y_n\Parrow 0$. 
  \end{enumerate}

\subsection{Configuration Random Graphs}\label{Configuration} 
\label{sec:1.2}
In their classic paper \cite{ErdoReny59}, Erd\"os and Renyi introduced the undirected model $G(N,M)$ that consists of   $N$ nodes and a random subset of exactly $M$ edges chosen uniformly from the collection of $\binom{N}{M}$ possible such edge subsets. This model can be regarded as the $M$th step of a random graph process that starts with $N$ nodes and no edges, and adds edges one at a time selected uniformly randomly from the set of available undirected edges. Gilbert's random graph model $G(N,p)$, which takes $N$ nodes and selects each possible edge independently with probability $p=z/(N-1)$, has mean degree $z$  and similar large $N$ asymptotics provided $M=zN/2$.   In fact, it was proved by \cite{Bollobas01} and \cite{MollReed95} that the undirected Erd\"os-Renyi graph  $G(N,zN/2)$ and $G(N,p_N)$ with probability $p_N=z/(N-1)$ both converge in probability to the same model as $N\to\infty$ for all $z\in\BBR_+$. Because of their popularity, the two models $G(N,p)\sim G(N,zN/2)$ have come to be known as ``the'' random graph\index{random graph}. Since the degree distribution of $G(N,p)$ is ${\rm Bin}(N-1, p)\sim_{N\to\infty} {\rm Pois}(z)$, this is also called the {\it Poisson graph} model\index{random graph!Poisson graph}. Both these constructions have obvious directed graph analogues.

The well known directed configuration multigraph model introduced by Bollobas \cite{Bollobas80} with general degree distribution $P=\{P_{jk}\}_{j,k=0,1,\dots}$ and size $N$ is constructed by the following random algorithm: \begin{enumerate}
  \item Draw a sequence of  $N$ node-type pairs $(j_1,k_1),\dots, (j_N,k_N)$ independently from $P$, and accept the draw if and only if it is feasible, i.e. $\sum_{n\in [N]}(j_n-k_n)=0$. Label the $n$th node with $k_n$ {\it out-stubs}\index{out-stub} (picture this as a half-edge with an out-arrow) and $j_n$ {\it in-stubs}\index{in-stub}.  
   \item While there remain available unpaired stubs, select (according to any rule, whether random or deterministic) any unpaired out-stub and pair it with an in-stub selected uniformly amongst unpaired in-stubs. Each resulting pair of stubs is a directed edge of the multigraph. 
\end{enumerate}

The algorithm leads to objects with self-loops and multiple edges, which are usually called multigraphs rather than graphs.  Only multigraphs that are free of self-loops and multiple edges, a condition called ``simple'', are considered to be graphs. For the most part, one does not care over much about the distinction, because the density of self-loops and multiple edges goes to zero as $N\to\infty$. In fact, Janson \cite{Janson09_simple} has proved in the undirected case that the probability for a multigraph to be simple is bounded away from zero for well-behaved sequences $(g_N)_{N>0}$ of size $N$ graphs with given $P$. 

Exact simulation of the adjacency matrix in the configuration model with general $P$ is problematic because the feasibility condition met in the first step occurs only with asymptotic frequency $\sim \frac{\sigma}{\sqrt{2\pi N}} $, which is vanishingly small for large graphs. For this reason, practical Monte Carlo implementations use some kind of rewiring or {\it clipping}\index{clipping}  to adjust each infeasible draw of node-type pairs. 

Because of the uniformity of the matching in step 2 of the above construction, the edge-type distribution of the resultant random graph is 
\be\label{independentedge} Q_{kj}=\frac{jkP^+_kP^-_j}{z^2}=Q^+_kQ^-_j\ee
which is called the {\it independent edge condition}\index{independent edge condition}. For many reasons, financial and otherwise, one is interested in the more general situation called {\it assortativity} when \eqref{independentedge} is not true. We will now show how such an extended class of assortative configuration graphs can be defined. The resultant class encompasses all reasonable type distributions $(P,Q)$ and has special properties that make it suitable for exact analytical results, including the possibility of a detailed percolation analysis. 

  \section{The ACG Construction}\label{finite_assortative_section}  
  \label{sec:2}
The assortative configuration (multi-)graph (ACG)  of size $N$ parametrized by the node-edge degree distribution pair $(P,Q)$ that satisfy the consistency conditions \eqref{Consistent1} is defined by the following random algorithm: \begin{enumerate}
  \item Draw a sequence of  $N$ node-type pairs $X= ((j_1,k_1),\dots, (j_N,k_N))$ independently from $P$, and accept the draw if and only if it is feasible, i.e. $\sum_{n\in [N]}\ j_n =\sum_{n\in [N]}\ k_n$, and this defines the number of edges $E$ that will result. Label the $n$th node with $k_n$ {\it out-stubs}\index{out-stub} (picture each out-stub as a half-edge with an out-arrow, labelled by its degree $k_n$) and $j_n$ {\it in-stubs}\index{in-stub}, labelled by their degree $j_n$.  Define the partial sums $u^-_j=\sum_n\One(j_n=j), u^+_k=\sum_n\One(k_n=k), u_{jk}=\sum_n\One(j_n=j,k_n=k)$, the number $ e_k^+=k u_k^+$ of {\it k-stubs} (out-stubs of degree $k$)  and the number of {\it j-stubs} (in-stubs of degree $j$), $ e_j^-=j u_j^-$.
   \item Conditioned on $X$, the result of Step 1, choose an arbitrary ordering $\ell^-$ and $\ell^+$ of the $E$ in-stubs and $E$ out-stubs. The matching sequence, or ``wiring'', $W$ of edges  is selected by choosing a pair of permutations $\sigma,\tilde\sigma\in S(E)$ of the set $[E]$. This determines the edge sequence $\ell=(\ell^-=\sigma(\ell), \ell^+=\tilde\sigma(\ell))$   labelled by  $\ell\in[E]$, to which is attached a probability weighting factor
   \be\label{AWweighting}\prod_{\ell\in[E]} Q_{k_{\sigma(\ell)}j_{\tilde\sigma(\ell)}}\ .
   \ee
\end{enumerate}

Given the wiring $W$ determined in Step 2, the number of type $(k,j)$ edges is 
\be\label{ekj}
e_{kj}=e_{kj}(W)=\sum_{\ell\in[E]}\One(k_{\tilde\sigma(\ell)}=k,\ j_{\sigma(\ell)}=j)\ .\ee
The collection $e=(e_{kj})$ of edge-type numbers are constrained by the $e^+_{k}, e^-_j$ that are determined by Step 1:
\be\label{e_constraints}
e^+_{k}=\sum_{j}e_{kj},\quad e^-_j=\sum_{k}e_{kj}, \quad E=\sum_{kj}e_{kj}\ .\ee

Intuitively, since Step 1 leads to a product probability measure subject to a single linear constraint that is true in expectation, one expects that it will lead to the independence of node degrees for large $N$, with the probability $P$. 
Similar logic suggests that since the matching weights in Step 2 define a product probability measure conditional on a set of linear constraints that are true in expectation, it should lead to edge degree independence in the large $N$ limit, with the limiting probabilities given by $Q$. However, the verification of these facts is not so easy, and their justification is the main object of this paper. First, certain combinatorial properties of the wiring algorithm of Step 2, conditioned on the node-type sequence $X$ resulting from Step 1 for a  finite $N$ will be derived. One result says that the probability of any wiring sequence $W=(\ell\in[E])$ in step 2 depends only on the set of quantities $(e_{kj})$ where for each $k,j$, $e_{kj}:=|\{\ell\in[E]|\ell\in\CE_{kj}\}|$. Another is that the conditional expectation of $e_{kj}/E$ is the exact edge-type probability for all edges in $W$. 

\begin{proposition}
\label{P1} Consider Step 2 of the ACG construction for finite $N$ with probabilities $P,Q$ conditioned on the $X=( j_i, k_i), i\in[N]$.
\begin{enumerate}
\item  The conditional probability of any wiring sequence $W=(\ell\in[E])$  is:
\beq\label{wiringprob}
\BBP[W|X]&=&C^{-1}\ \prod_{kj} (Q_{kj})^{e_{kj}(W)}\ ,\\
\label{Cformula}C&=&C(e^-,e^+)=E! \sum_{e}\prod_{kj} \frac{(Q_{kj})^{e_{kj}}}{e_{kj}!}\prod_j\left(e_j^-!\right)\prod_k\left(e_k^+!\right)\ ,
\eeq
where the sum in \eqref{Cformula} is over collections $e=(e_{kj})$ satisfying the constraints \eqref{e_constraints}. \item  The conditional probability $p$ of any edge of the wiring sequence $W=(\ell\in[E])$ having type $k, j$ is  
\be\label{E_prob}p=\BBE[e_{kj}|X]/E\ .\ee
\end{enumerate} 
\end{proposition}

\bigskip\noindent{\bf Proof of Proposition \ref{P1}: \ } The denominator of \eqref{wiringprob} is $C=\sum_{\sigma,\tilde\sigma\in S(E)}\prod_{l\in[E]}Q_{k_{\sigma(\ell)}j_{\tilde\sigma(\ell)}}$, from which \eqref{Cformula} can be verified by induction on $E$.
Assuming \eqref{Cformula} is true for $E-1$, one can verify the inductive step for $E$:
\beqq  C&=&\sum_{\tilde k,\tilde j} \sum_{\sigma,\tilde\sigma\in S(E)}\One(k_{\sigma(E)}=\tilde k,j_{\tilde\sigma(E)}=\tilde j)\prod_{l\in[E]}\ Q_{k_{\sigma(\ell)}j_{\tilde\sigma(\ell)}}\\
&=&\sum_{\tilde k,\tilde j} e^+_{\tilde k}e^-_{\tilde j}\ Q_{\tilde k\tilde j}\ \sum_{\sigma',\tilde\sigma'\in S(E-1)}\prod_{l\in[E-1]}\ Q_{k_{\sigma'(\ell)}j_{\tilde\sigma'(\ell)}}\\
&=&\sum_{\tilde k,\tilde j} e^+_{\tilde k}e^-_{\tilde j}\ Q_{\tilde k\tilde j}\ (E-1)!\  \sum_{e'}\prod_{kj} \frac{(Q_{kj})^{e'_{kj}}}{e'_{kj}!}\prod_j\left(e_j^{'-}!\right)\ \prod_k\left(e_k^{'+}!\right)\ . \eeqq
Here, $e'_{kj}=e_{kj}-\One(k=\tilde k,j=\tilde j),\ e_j^{'-}=e_j^--\One(j=\tilde j), \ e_k^{'+}=e_k^+-\One(k=\tilde k)$. After noting cancellations that occur in the last formula, and re-indexing the collection $e'$ one finds
\beqq  C&=&\sum_{\tilde k,\tilde j}  \sum_{e'}e_{\tilde k\tilde j} \  (E-1)!\prod_{kj} \frac{(Q_{kj})^{e_{kj}}}{e_{kj}!}\prod_j\left(e_j^-!\right)\ \prod_k\left(e_k^+!\right)\ \\
&=& \sum_{e}\left(\sum_{\tilde k,\tilde j} e_{\tilde k\tilde j}\right)\  (E-1)!\ \prod_{kj} \frac{(Q_{kj})^{e_{kj}}}{e_{kj}!}\prod_j\left(e_j^-!\right)\ \prod_k\left(e_k^+!\right)\ 
\\ &=&E!\sum_{e}\prod_{kj} \frac{(Q_{kj})^{e_{kj}}}{e_{kj}!}\prod_j\left(e_j^-!\right)\ \prod_k\left(e_k^+!\right)\ \eeqq
which is the desired result. 

Because of the edge-permutation symmetry, it is enough to prove \eqref{E_prob} for the last edge. For this, one can follow the same logic and steps as in Part 1 to find:
\beqq p&=& \frac1{C(e^-,e^+)}
 \sum_{\sigma,\tilde\sigma\in S(E)}\One(k_{\sigma(E)}=k,j_{\tilde\sigma(E)}=j)\ \prod_{l\in[E]}\ Q_{k_{\sigma(\ell)}j_{\tilde\sigma(\ell)}}\\&=&\frac{E!}{C(e^-,e^+)}  \sum_{e} \frac{e_{kj}}{E}\ \prod_{k'j'} \frac{(Q_{k'j'})^{e_{k'j'}}}{e_{k'j'}!}\prod_{j'}\left(e_{j'}^{-}!\right)\ \prod_{k'}\left(e_{k'}^{+}!\right)=
 \BBE[e_{kj}|X]/E\ .\eeqq
  \qquad

\qedsymbol

\bigskip 

An easy consequence of the above proof is that the number of wirings $W$ consistent with a collection $e=(e_{kj})$ is given by 
\be\label{Nwiring} |\{W:e(W)=e\}|=\frac{E!\left(\prod_je^-_j!\right)\left(\prod_k e^+_k!\right)}{\prod_{kj}e_{kj}!}\ .\ee

Because of the permutation symmetries of the construction, a host of more complex combinatorial identities hold for this model. The most important is that Part 2 of the Proposition can be extended inductively to determine the joint edge distribution for the first $M$ edges conditioned on $X$. To see how this goes, define two sequences $e^-_j(m), e^+_k(m)$ for $0\le m\le M$ to be the number of  $j$-stubs and $k$-stubs available after $m$ wiring steps.
\begin{proposition}
\label{P2} Consider Step 2 of the ACG construction for finite $N$ with probabilities $P,Q$ conditioned on $X$ from Step 1.
The conditional probability $p$ of the first $M$ edges of the wiring sequence $W=(\ell\in[E])$ having types $(k_i, j_i)_{i\in[M]}$ is  
\be\label{EM_prob}\BBP[(k_i, j_i)_{i\in[M]}|X]=\frac{(E-M)!}{E!}\prod_{i\in[M]}\BBE[e_{k_ij_i}|e^-(i-1), e^+(i-1)]\ .\ee
\end{proposition}

\bigskip\noindent{\bf Proof of Proposition \ref{P2}: \ } Note that Part 2 of Proposition \ref{P1} gives the correct result when $M=1$. For any $m$, an extension of the argument that proves Part 2 of Proposition \ref{P1} also shows that
\beq\BBP[(k_i, j_i)_{i\in[m]}|X]&=& \frac1{C(e^-(0),e^+(0))}
 \sum_{\sigma,\tilde\sigma\in S(E)}\prod_{\ell=1}^m\One(k_{\sigma(\ell)}=k_\ell,j_{\tilde\sigma(\ell)}=j_\ell)\ \prod_{l\in[E]}\ Q_{k_{\sigma(\ell)}j_{\tilde\sigma(\ell)}}\nonumber\\
 &&\hspace{-1.5in} =\ \frac{1}{C(e^-(0),e^+(0))}  \prod_{\ell=1}^m\left[ e^-_{j_\ell}(\ell-1)e^+_{k_\ell}(\ell-1)\ Q_{k_\ell j_\ell}\right] \ \sum_{\sigma',\tilde\sigma'\in S(E-m)}\prod_{\ell=m+1}^{E}\ Q_{k_{\sigma'(\ell)}j_{\tilde\sigma'(\ell)}}
 \label{EM_prob2}\eeq
Now assume inductively that the result  \eqref{EM_prob} is true for $M-1$ and compute \eqref{EM_prob} for $M$:
\[\BBP[(k_i, j_i)_{i\in[M]}|X]= \frac{\BBP[(k_i, j_i)_{i\in[M]}|X]}{\BBP[(k_i, j_i)_{i\in[M-1]}|X]}\times \frac{(E-M-1)!}{E!}\prod_{i\in[M-1]}\BBE[e_{k_ij_i}|e^-(i-1), e^+(i-1)]\ .
\]
The ratio in the first factor can be treated using \eqref{EM_prob2}, and the resulting cancellations lead to the formula
\beqq
\BBP[(k_i, j_i)_{i\in[M]}|X]&&\\&&\hspace{-1.5in}=\ \frac{\left[ e^-_{j_M}(M-1)e^+_{k_M}(M-1)\ Q_{k_M j_M}\right] \ \sum_{\sigma',\tilde\sigma'\in S(E-M)}\prod_{\ell=M+1}^{E}\ Q_{k_{\sigma'(\ell)}j_{\tilde\sigma'(\ell)}}}{ \ \sum_{\sigma',\tilde\sigma'\in S(E-M+1)}\prod_{\ell=M}^{E}\ Q_{k_{\sigma'(\ell)}j_{\tilde\sigma'(\ell)}}}\\
&&\hspace{-1.5in}\times\ \frac{(E-M+1)!}{E!}\prod_{i\in[M-1]}\BBE[e_{k_ij_i}|e^-(i-1), e^+(i-1)]
\eeqq
The desired result follows because Part 2 of Proposition \ref{P1} can be applied to show
\beqq \frac{\left[ e^-_{j_M}(M-1)e^+_{k_M}(M-1)\ Q_{k_M j_M}\right] \ \sum_{\sigma',\tilde\sigma'\in S(E-M)}\prod_{\ell=M+1}^{E}\ Q_{k_{\sigma'(\ell)}j_{\tilde\sigma'(\ell)}}}{ \ \sum_{\sigma',\tilde\sigma'\in S(E-M+1)}\prod_{\ell=M}^{E}\ Q_{k_{\sigma'(\ell)}j_{\tilde\sigma'(\ell)}}}\\
&&\hspace{-2.7in}=\frac1{E-M+1}\BBE[e_{k_Mj_M}|e^-(M-1), e^+(M-1)]\ .
\eeqq
\qedsymbol

\bigskip 

\section{Asymptotic Analysis}

It is quite easy to prove that the empirical node-type distributions $(u_{jk},u^-_j,u^+_k)$  resulting from Step 1 of the ACG algorithm satisfy a law of large numbers: 
\be\label{LLN_P}
N^{-1}u_{jk}\Pequals P_{jk},\quad N^{-1}u^-_{j}\Pequals P^-_{j},\quad N^{-1}u^+_{k}\Pequals P^+_{k},
\ee
 as $N\to\infty$. In this section, we focus on the new and more difficult problem to determine the asymptotic law of the empirical edge-type distribution, conditioned on the node-type sequence $X$. To keep the discussion as clear as possible, we confine the analysis to the case the distributions $P$ and $Q$ have support on the finite set $(j,k)\in\{0,1,\dots, K\}^2$.
\ 
One can see from Proposition \ref{P2} that the probability distribution of the first $M$ edge types will be given  asymptotically by $\prod_{i\in[M]}Q_{k_ij_i}$ provided our intuition is correct that $\BBE[E^{-1}e_{kj}]\Pequals Q_{kj}(1+o(1))$ asymptotically for large $N$. To validate this intuition, it turns out one can apply the Laplace asymptotic method to the cumulant generating function for the empirical edge-type random variables $e_{kj}$, conditioned on any feasible collection of $(e^+_k,e^-_j)$ with total number $E=\sum_ke^+_k=\sum_j e^-_j$:
\beq\label{cumfn}
F(v;e^-,e^+)&:=&\log\BBE[e^{\sum_{kj}v_{kj} e_{kj}}|e^-,e^+],  \quad \forall v=(v_{kj})\\
&=&\log \frac{\sum_{e}\prod_{kj} \frac{(Q_{kj}e^{v_{kj}})^{e_{kj}}}{e_{kj}!}\prod_j\left(e_j^-!\right)\prod_k\left(e_k^+!\right)}
{\sum_{e}\prod_{kj} \frac{(Q_{kj})^{e_{kj}}}{e_{kj}!}\prod_j\left(e_j^-!\right)\prod_k\left(e_k^+!\right)}\ ,
\eeq
The constraints on $(e_{kj})$ can be introduced by auxiliary integrations over variables $u^-_j,u^+_k$ of the form
\[ 2\pi\One(\sum_j e_{kj}=e_k^+)=\int^{2\pi}_0 d u^+_k \ e^{i\tilde u_k(\sum_j e_{kj}-e_k^+)}\ .
\]
This substitution leads to closed formulas for the sums over $e_{kj}$,  and the expression for $e^F$:
\be
\label{steep}\frac{\int_I d^{2K}u \exp[H(v,-iu;e)]}{\int_I d^{2K}u\ \exp [H(0,-iu;e)]}
\ee
where 
\be\label{Hfunction} H(v,\alpha;e)=\sum_{kj}e^{(\alpha^-_j+\alpha^+_k)}e^{v_{kj}} Q_{kj} -(\sum_j \alpha^-_je^-_j+\sum_k\alpha^+_k e^+_k)=\sum_{kj}e^{\alpha\cdot\delta_{jk}}e^{v_{kj}} Q_{kj} -\alpha\cdot e\ .
\ee
The integration in \eqref{steep}  is over the set $I:=[0,2\pi]^{2K}$.

Here a ``double vector'' notation has been introduced for  $u=(u^-,u^+), e=(e^-,e^+), \alpha=(\alpha^-,\alpha^+)$ where $u^-, u^+\in\BBC^K$ etc. and where $K$ is the number of possible in and out degrees (which one may want to take to be infinite). Define double vectors $\One^-=(1,1,\dots\ 1;0,\dots,\ 0), \One^+=(0,\dots,\ 0;1,\dots,\ 1), \One=\One^-+\One^+,\tilde\One=\One^--\One^+$. For any pair $(j,k)\in [K]^2$, let $\delta^-_j$ be the double vector with a $1$ in the $j$th place and zeros elsewhere, let  $\delta^+_k$ be the double vector with a $1$ in the $K+k$th place and zeros elsewhere and $\delta_{jk}=\delta^-_j+\delta^+_k$. Using the natural inner product for double vectors  $\alpha\cdot e:=\sum_j\alpha^-_je^-_j+\sum_k \alpha^+_ke^+_k$, etc., the feasibility condition on stubs can be written $e\cdot\tilde\One=0$.

The main aim of the paper is to prove a conditional law of large numbers for $E^{-1}e_{jk}$ as $E\to\infty$, conditioned on $e=(e^-,e^+)$ satisfying $e\cdot\tilde\One=0$. By explicit differentiation of the cumulant generating function, and some  further manipulation, one finds that
\beq
\BBE[e_{kj}|e]&=&\frac{\pa F}{\pa v_{kj}}\Big|_{v=0}=Q_{kj}\frac{\int_I d^{2K}u\  \exp[H(0,-iu;e-\delta_{jk})]}{\int_I d^{2K}u \ \exp [H(0,-iu;e)]}\label{firstmom}\\
\var[e_{kj}|e]&=&\frac{\pa^2 F}{\pa v_{kj}^2}\Big|_{v=0}=Q_{kj}\frac{\int_I d^{2K}u\  \exp[H(0,-iu;e-\delta_{jk})]}{\int_I d^{2K}u\ \exp [H(0,-iu;e)]}
\label{secondmom}
\\
&&\hspace{-1.05in}+\left(Q_{kj}\right)^2\left[\frac{\int_I d^{2K}u\  \exp[H(0,-iu;e-2\delta_{jk})]}{\int_I d^{2K}u \ \exp [H(0,-iu;e)]}-\left(\frac{\int_I d^{2K}u\  \exp[H(0,-iu;e-\delta_{jk})]}{\int_I d^{2K}u\ \exp [H(0,-iu;e)]}\right)^2\right]\ .\nonumber
\eeq

Since our present aim is to understand \eqref{firstmom} and \eqref{secondmom}, we henceforth set $v=0$ in the $H$-function. The  $H$ function defined by \eqref{Hfunction} with $v=0$ has special combinatorial features:
\begin{lemma}\label{Lem1} For all $e\in\BBZ_+^{2K}$ satisfying $e\cdot\tilde\One=0$, the function $H=H(\alpha; e)$ satisfies the following properties:
\begin{enumerate}
\item $H$ is convex for $\alpha\in\BBR^{2K}$ and entire analytic for $\alpha\in\BBC^{2K}$;
\item $H$ is periodic: $H(\alpha+2\pi i \eta;e)=H(\alpha; e)$ for all $\eta\in\BBZ^{2K}$.
\item   For any $\lambda\in \BBC$, $H(\alpha+\lambda\tilde\One;e)=H(\alpha; e)$  ;
  \item For any $\lambda>0$, $H(\alpha;\lambda e)=\lambda H(\alpha-\frac{\log\lambda}{2}\One;e)-\lambda\ \log\lambda \One\cdot e
 $.
\item The $m$th partial derivative of $H$ with respect to $\alpha$ is given by
\be \nabla^m H(\alpha;e) =\left\{\begin{array}{ll }
     \sum_{jk} \delta_{jk} e^{\alpha\cdot\delta_{jk}}\ Q_{kj}-e,&  m=1;    \\
       \sum_{jk} (\delta_{jk})^{\bigotimes m} e^{\alpha\cdot\delta_{jk}}\ Q_{kj},& m=2,3,\dots 
\end{array}
\right.\ee Here $(\delta_{jk})^{\bigotimes m}$ denotes the $m$th tensor power of the double vector $\delta_{jk}$. 
\end{enumerate}
\end{lemma}

The Laplace asymptotic method (or saddlepoint method), reviewed for example in \cite{Hurd10}, involves shifting the $u$ integration in \eqref{steep}  into the complex by an imaginary vector. The Cauchy Theorem, combined with the periodicity of the integrand in $u$, will ensure the value of the integral is unchanged under the shift.  The desired shift is determined by the $e$-dependent critical points $\alpha^*$ of $H$  which by Part 5 of Lemma  \ref{Lem1} are solutions of
\be\label{criticalpt} \sum_{jk} \delta_{jk} e^{\alphabf\cdot\delta_{jk}}\ Q_{kj} = e\ . 
\ee
In view of Parts 1 and 3 of the Lemma, for each $e\in\BBZ^{2K}$ there is a unique critical point $\alpha^*(e)$ such that $\tilde\One\cdot\alpha^*(e)=0$. The imaginary shift of the $u$-integration is implemented by writing $u=i\alpha^*(e)+\zeta$ where now $\zeta$ is integrated over $I$.

To unravel the $E$ dependence, one uses rescaled variables $x=E^{-1} e$ that lie on the plane $\tilde\One\cdot x=1$ and by Part 4 of the Lemma with $\lambda=E^{-1}$ one has that $\alpha^*(e)=\alpha^*(x)+\frac{\log E}{2}\ \One$.  Now one can use the third order Taylor expansion with remainder to write
\beq H(0,\alpha^*(e)-i\zeta;e)&=& EH(0,\alpha^*(x)-i\zeta;x)-E\log E (\One\cdot x)\label{taylor} \\
&&\hspace{-1.2in}=\ -\ E\log E+E\left[H(0,\alpha^*(x);x)-\frac12\zeta^{\bigotimes 2}\cdot\nabla^2H+i\frac16\zeta^{\bigotimes 3}\cdot\nabla^3H\right]+EO(|\zeta|^4)\nonumber
\eeq
where $\nabla^2H, \nabla^3H$ are evaluated at $\alpha^*(x)$ and the square-bracketed quantities are all $E$ independent. From \eqref{Hfunction} one can observe directly that $|e^H|$ has a unique maximum on the domain of integration at $\zeta=0$:
\be\label{maxonI}
\max_{\zeta\in I}|e^{H(\alpha^*(e)-i\zeta;e)}|=e^{H(\alpha^*(e);e)}\ .
\ee
The uniqueness of the maximum is essential to validate the following Laplace asymptotic analysis, and leads to the main result of the paper:
\begin{theorem}
\label{theorem}
 For any double vector $x^*\in(0,1)^{2K}\cap \tilde\One^\perp$, let $e(E)=Ex(E)$ be a sequence in $\BBZ_+^{2K}\cap \tilde\One^\perp$ such that 
 \be\lim_{E\to\infty} x(E)=x^*\ .
 \ee
 Then asymptotically as $E\to\infty$, 
  \beq \CI(E)&=&\int_I d^{2K}u\  \exp[H(-iu;e(E))]\\
  & =& (2\pi)^{K+1/2}E^{1/2-K}e^{-E\log E+ EH(\alpha^*(x^*);x^*)}\left[{\det}_0\nabla^2H\right]^{-1/2}\left[1+O(E^{-1})\right]\ .\nonumber
  \eeq
 Here ${\det}_0\nabla^2H$ represents the determinant of the matrix projection onto $\tilde\One^\perp$, the subspace orthogonal to $\tilde\One$, of $\nabla^2H$ evaluated at the critical point $\alpha^*(x^*)$. 
 
\end{theorem}

When applied to \eqref{firstmom} and \eqref{secondmom} this Theorem is powerful enough to yield the desired results on the edge-type distribution in the ACG model for fixed $e=(e^-,e^+)=Ex$ for large $E$.
 \begin{corollary}\label{corollary}
 Consider the ACG model with $(P,Q)$ supported on $\{0,1,\dots, K\}^2$.\begin{enumerate}
  \item Conditioned on $X$,
  \[ E^{-1}e_{kj}\Pequals Q_{kj}e^{1-H(0,\alpha^*(x);x)-\alpha^*(x)\cdot\delta_{jk}}[1+O(E^{-1/2})]\]
where $x=E^{-1}e$ and $e=(e^-(X),e^+(X))$.
\item Unconditionally, 
\[ E^{-1}e_{kj}\Pequals Q_{kj}[1+O(N^{-1/2})]\ .\]
\end{enumerate}\end{corollary}

Combining this Law of Large Numbers result with the easier result for the empirical node-type distribution confirms that the large $N$ asymptotics of the empirical node- and edge-type distributions agree with the target $(P,Q)$ distributions. 

\bigskip\noindent{\bf Proof of Corollary \ref{corollary}: \ } By applying Part 4 of Lemma \ref{Lem1} and the Theorem to \eqref{firstmom} one finds that
\beqq \BBE[e_{kj}|e]&=&Q_{kj}\frac{\int_I d^{2K}u\  \exp[H(-iu;e-\delta_{jk})]}{\int_I d^{2K}u \ \exp [H(-iu;e)]}\\
&&\hspace{-.9in}=\ Q_{kj}\ \exp[-(E-1)\log (E-1)+ E\log E+ (E-1)H(\alpha^*(x');x')-EH(\alpha^*(x);x)]\\
&&\times\ \left[\frac{{\det}_0\nabla^2H(\alpha^*(x))}{{\det}_0\nabla^2H(\alpha^*(x'))}\right]^{1/2}\left[1+O(E^{-1})\right]\eeqq
where $x=E^{-1}e$ and $x'=(E-1)^{-1}(e-\delta_{jk})$ are such that $\Delta x=x'-x=O(E^{-1})$. Now,
one can show that if $x, x'$ lie on the plane $\tilde\One\cdot x=1$, and $\Delta x=x'-x$ is $O(E^{-1})$ then
 \be\label{DeltaH} H(\alpha^*(x');x')-H(\alpha^*(x);x)=\alpha^*(x)\cdot\Delta x +O(|\Delta x|^2)\ .
 \ee
It is also true that $\Delta\alpha^*=\alpha^*(x')-\alpha^*(x)=O(|\Delta x|)$ and satisfies 
 \be \Delta\alpha^*\cdot x=O(|\Delta x|^2)\ .\ee
Since ${\det}_0\nabla^2H(\alpha)$ is analytic in $\alpha$ with $O(1)$ derivatives, and $\Delta\alpha^*=O(|\Delta x|)$ \[ \left[\frac{{\det}_0\nabla^2H(\alpha^*(x))}{{\det}_0\nabla^2H(\alpha^*(x'))}\right]^{1/2}=\left[1+O(E^{-1})\right]\ .
\]
Also, 
\[(E-1)H(\alpha^*(x');x')-(E-1)H(0,\alpha^*(x);x)=-\alpha^*(x)\cdot\delta_{jk} +O(|\Delta x|)\]
and $E\log E-(E-1)\log (E-1) ]\sim \log E +1 +O(E^{-1})$, from which one concludes
\be\BBE[e_{kj}|e]=Q_{kj}\ E\ \exp[1-H(\alpha^*(x);x)-\alpha^*(x)\cdot\delta_{jk}]\left[1+O(E^{-1})\right]\ .\ee

The conclusion of the Part 1 of the Corollary now follows from the Chebyshev inequality if one shows that \eqref{secondmom} is $O(E)$. Since the first term of \eqref{secondmom} equals $\BBE[e_{kj}|e]$, which is $O(E)$, it is only necessary to show that the $O(E^2)$ parts of the second term cancel. Each ratio in the second term can be analyzed exactly as above, leading to  
\beqq  (Q_{kj})^2
E\ \exp[1-H(\alpha^*(x);x)-\alpha^*(x)\cdot\delta_{jk}]&&\\&&\hspace{-3.25in}\times \ \left((E-1) \exp[1-H(\alpha^*(x');x')-\alpha^*(x')\cdot\delta_{jk}]-E\ \exp[1-H(\alpha^*(x);x)-\alpha^*(x)\cdot\delta_{jk}]\right)\\
&&\hspace{-2.5in}\times \ \left[1+O(E^{-1})\right]\\
=\ \left[Q_{kj}E \exp[1-H(\alpha^*(x);x)-\alpha^*(x)\cdot\delta_{jk}]\right]^2&&\\&&\hspace{-2.3in}\times \ \left(\exp[H(\alpha^*(x);x)-H(\alpha^*(x');x')-\Delta\alpha^*(x')\cdot\delta_{jk}]\right)\left[1+O(E^{-1})\right]\\
=\ \left[Q_{kj}E \exp[1-H(\alpha^*(x);x)-\alpha^*(x)\cdot\delta_{jk}]\right]^2&&\\&&\hspace{-2.3in}\times \ \left(\exp[-\alpha^*(x)\cdot\Delta x -\Delta\alpha^*(x')\cdot\delta_{jk}]-1\right)\left[1+O(E^{-1})\right]\ =\ O(E)
\eeqq
where one uses \eqref{DeltaH} again in the second last equality.

To prove Part 2, it is sufficient to note that $E^{-1}(e^-(X),e^+(X))=(Q^-,Q^+)[1+O(N^{-1/2})]$ and that $\alpha^*(Q^-,Q^+)=0, H(\alpha^*(Q^-,Q^+);Q^-,Q^+)=1$. 

\qedsymbol

\bigskip
\noindent{\bf Proof of Theorem \ref{theorem}: \ } For each $E$, since the integrand of  $\CI(E)$ is entire analytic and periodic, its integral is unchanged under a purely imaginary shift of the contour. Also, since by Part 3 of Lemma \ref{Lem1} the integrand is constant in directions parallel to $\tilde\One$, the integrand can be reduced to the set $I\cap\tilde\One^\perp$. Thus, using \eqref{taylor} for $e=e(E)$ and $x=x(E)$, $\CI(E)$ can be written
\beqq
\CI(E)&=&2\pi\int_{I\cap\tilde\One^\perp} d^{2K-1}\zeta\  \exp[H(\alpha^*(e)-i\zeta;e)]\\
&=& 2\pi\int_{I\cap\tilde\One^\perp} d^{2K-1}\zeta\\
&&\hspace{-.8in}\times \   \exp\left[-\ E\log E+E\left(H(\alpha^*(x);x)-\frac12\zeta^{\bigotimes 2}\cdot\nabla^2H+i\frac16\zeta^{\bigotimes 3}\cdot\nabla^3H+O(|\zeta|^4)\right)\right]\ .
\eeqq
In rescaled variables $\tilde\zeta=E^{1/2}\zeta$ this becomes
\[
\CI(E)=2\pi E^{1/2-K}\exp\left[-\ E\log E+EH(\alpha^*(x);x)\right]\times \tilde \CI(E)\]
where
\beqq
\tilde \CI(E)&:=&\int_{E^{1/2}I\cap\tilde\One^\perp} d^{2K-1}\tilde\zeta\ \exp[ H(\alpha^*(x)-iN^{-1/2}\tilde\zeta;x)- H(\alpha^*(x);x)]\\
&&\hspace{-.5in}=\ \int_{E^{1/2}I\cap\tilde\One^\perp} d^{2K-1}\tilde\zeta\ \exp[-\frac12\tilde\zeta^{\bigotimes 2}\cdot\nabla^2H]\left(1+i\frac{E^{-1/2}}{6} \tilde\zeta^{\bigotimes 3}\cdot\nabla^3H+O(E^{-1})\right)\ .
\eeqq
In this last integral the $O(E^{-1/2})$ term is odd in $\tilde\zeta$ and makes no contribution. Now,  
\beqq |\exp[H(\alpha^*(x)-iN^{-1/2}\tilde\zeta;x)- H(\alpha^*(x);x)]|&&\\
&&\hspace{-2in}=\ \exp\left[\sum_{kj}e^{\alpha^*(x)\cdot\delta_{jk}}\ (\cos(N^{-1/2}\tilde\zeta)-1)\cdot\delta_{jk})\ Q_{kj}\right]
\eeqq clearly has a unique maximum at $\tilde\zeta=0$. Therefore, a standard version of the Laplace method such as that found in \cite{Erdelyi1956} is sufficient to imply that as $E\to\infty$, 
\be \tilde \CI(E)= (2\pi)^{K-1/2}\left[{\det}_0\nabla^2H\right]^{-1/2}\left[1+O(E^{-1})\right]\ee
where $\nabla^2H$ is evaluated at $\alpha^*(x^*)$.

\qedsymbol

\section{Locally Tree-like Property} \label{Sec:4}

To understand percolation theory on random graphs, or to derive a rigorous treatment of cascade mappings on random financial networks, it turns out to be important that the underlying random graph model have a property sometimes called ``locally tree-like''. In this section, the local tree-like property of the ACG model will be characterized as a particular large $N$ property of  the probability distributions associated with {\it configurations}, that is, finite connected subgraphs $g$ of the skeleton labelled by their degree types.

First consider what it means in the $(P,Q)$ ACG model with size $N$ to draw a random configuration $g$ consisting of a pair of vertices $v_1,v_2$ joined by a link, that is, $v_2\in\CN^-_{v_1}$. In view of the permutation symmetry of the ACG algorithm, the random link can without loss of generality be taken to be the first link $W(1)$ of the wiring sequence $W$. Following the ACG algorithm, Step 1 constructs a feasible node degree sequence $X=( j_i, k_i), i\in[N]$ on nodes labelled by $v_i=i$ and conditioned on $X$, Step 2 constructs a random $Q$-wiring sequence $W=\left(\ell=(v^+_{\ell},v^-_{\ell})\right)_{\ell\in[E]}$ with $E=\sum_i 
 k_i=\sum_i  j_i$ edges. By an abuse of notation, we label their edge degrees by $k_{\ell}=k_{v^+_{\ell}},  j_{\ell}=j_{v^-_{\ell}}$ for $\ell\in[E]$ . 
The configuration event in question, namely that the first link in the wiring sequence $W$ attaches to nodes of the required degrees $(j_1,k_1), (j_2,k_2)$, has probability $p=\BBP[v_i\in\CN_{j_i,k_i}, i=1,2| v_2\in\CN^-_{v_1}]$. To compute this, note that  the fraction  $j_1u_{j_1k_1}/e^-_{j_1}$ of available $j_1$-stubs come from a $j_1k_1$ node and the fraction $k_2u_{j_2k_2}/e^+_{k_2}$ available $k_2$-stubs come from a $j_2k_2$ node. Combining this fact with Part 2 of Proposition 1, equation \ref{E_prob} implies the configuration probability conditioned on $X$ is 
exactly \be
p=j_{1}u_{j_{1}k_{1}}k_2u_{j_2k_2}\ \frac{\BBE[e_{k_2j_{1}}|e^-,e^+]}{Ee^+_{k_2}e^-_{j_{1}}}\ .
\ee
By the Corollary:
\be\label{one_animal} p\ \Pequals\  \frac{j_1k_2P_{j_1k_1}P_{j_2k_2} Q_{k_2j_1}}{z^2Q^+_{k_2}Q^-_{j_1}}[1+O(N^{-1/2})]\ .
\ee

This argument justifies the following informal computation of the correct asymptotic expression for $p$ by successive conditioning: 
\beq p&=&\BBP[v_i\in\CN_{j_ik_i}, i=1,2\left| v_2\in\CN^-_{v_1}]\right.\\
&=&\BBP[v_1\in\CN_{j_1k_1}\left| v_2\in\CN^-_{v_1}\cap\CN_{j_2k_2}]\right.\BBP[v_2\in\CN_{j_2k_2}\left| v_2\in\CN^-_{v_1}]\right.\\
&=&P_{k_1|j_1}Q_{j_1|k_2}P_{j_2|k_2}Q^+_{k_2}=\frac{P_{j_1k_1}P_{j_2k_2} Q_{k_2j_1}}{P^+_{k_2}P^-_{j_1}}
\eeq
where we introduce conditional degree probabilities $P_{k|j}=P_{jk}/P^-_j$ etc. 

Occasionally in the above matching algorithm, the first edge forms a self-loop, i.e. $v_1=v_2$. The probability of this event, jointly with fixing the degree of $v_1$, can be computed exactly for finite $N$ as follows:
\[ \tilde p:=\BBE[v_1=v_2, v_1\in\CN_{jk}|v_2\in\CN^-_{v_1}|X]=\left(\frac{jku_{jk}}{e^-_je^+_k}\right)\frac{\BBE[e_{kj}|X]}{E}\ .\]
As $N\to\infty$ this goes to zero, while $N\tilde p$ approaches a finite value:
\be\label{self_loop} N\tilde p\Parrow \frac{jkP_{jk}Q_{kj}}{z^2Q^+_kQ^-_j}
\ee
which says that the relative fraction of edges being self loops is the asymptotically small $\sum_{jk}\frac{jkP_{jk}Q_{kj}}{Nz^2Q^+_kQ^-_j}$. In fact, following results of \cite{Janson09_simple} and others on the undirected configuration model, one expects that the total number of self loops in the multigraph converges in probability to a Poisson random variable with finite parameter 
\be \lambda = \sum_{jk}\frac{jkP_{jk}Q_{kj}}{z^2Q^+_kQ^-_j}\ .
\ee

\subsection{General Configurations}
A general {\it configuration} is a connected subgraph $h$ of an ACG graph $(\CN,\CE)$ with $L$ ordered edges and with each node labelled by its degree type. It results from a growth process that starts from a fixed node $w_0$ called the root and at step $\ell\le L$ adds one edge $\ell$ that connects a node $w_\ell$ to a specific existing node $w'_\ell$. The following is a precise definition:
\begin{definition} A configuration rooted to a node $w_0$ with degree $(j,k):=(j_0,k_0)$ is a connected subgraph $h$ consisting of a sequence  of $L$ edges that connect nodes $(w_\ell)_{\ell\in[L]}$ of types
$(j_\ell,k_\ell)$, subject to the following condition:
For each $\ell\ge 1$, $w_\ell$ is connected by the edge labelled with $\ell$ to a node $w'_\ell\in\{w_j\}_{j\in\{0\}\cup[\ell-1]}$ by either an in-edge (that points into $w'_\ell$) $(w_\ell,w'_\ell)$ or an out-edge $(w'_\ell,w_\ell)$. 

\end{definition}

A random realization of the configuration results when the construction of the size $N$ ACG graph $(\CN,\CE)$ is conditioned on $X$ arising from Step 1 and the first $L$ edges of the wiring sequence of Step 2.  
The problem is to compute the probability of the node degree sequence $(j_\ell,k_\ell)_{\ell\in[L]}$ conditioned on $X$, the graph $h$ and the root degree $(j,k)$, that is
\be\label{configprob}
p=\BBP[w_\ell\in\CN_{j_\ell,k_\ell}, \ell\in[L]|w_0\in\CN_{jk},h,X] \ .
\ee

Note that there is no condition that the node $w_\ell$ at step $\ell$ is distinct from the earlier nodes $w_{\ell'}, \ell'\in{\{0\}\cup [\ell-1]}$. With high probability each $w_\ell$  will be new, and the resultant subgraph $h$ will be a tree with $L$ distinct added nodes (not including the root) and  $L$ edges. With small probability one or more of the $w_\ell$ will be preexisting, i.e. equal to $w_{\ell'}$ for some $\ell'\in{\{0\}\cup [\ell-1]}$: in this case the subgraph $h$ will have $M<L$  added nodes, will have cycles and not be a tree.

The following sequences of numbers are determined given $X$ and $h$: 
\begin{itemize}
\item $e_{j,k}(\ell)$ is the number of available $j$-stubs connected to $(j,k)$ nodes after $\ell$ wiring steps;
\item $e_{k,j}(\ell)$ is the number of available $k$-stubs connected to $(j,k)$ nodes after $\ell$ wiring steps.
\item $e^-_{j}(\ell):=\sum_ke_{j,k}(\ell)$ and $e^+_{k}(\ell):=\sum_je_{k,j}(\ell)$ are the number of available $j$-stubs and $k$-stubs respectively after $\ell$ wiring steps. 
\end{itemize}
Note that $e_{j,k}(0)=ju_{jk}$ and $e_{k,j}(0)=ku_{jk}$, and both decrease by at most $1$ at each step.

The analysis of configuration probabilities that follows is inductive on the step $\ell$.

\begin{theorem}\label{animal_theorem} Consider the  ACG sequence with $(P,Q)$  supported on $\{0,1,\dots, K\}^2$. Let $h$ be any fixed finite configuration rooted to $w_0\in\CN_{jk}$, with $M$ added nodes and $L\ge M$ edges, labelled by the  node-type sequence $( j_m, k_m)_{m\in[M]}$. Then,  as $N\to\infty$, the joint probability conditioned on $X$, 
\[p=\BBP[w_m\in\CN_{j_mk_m}, m\in[M]|w_0\in\CN_{jk},h,X]\ ,\]
is given by 
\be
 \prod_{m\in[M],\mbox{\ out-edge} }P_{k_m|j_{m}}Q_{j_m|k_{m'}} \prod_{m\in[M],\mbox{\ in-edge} }P_{j_m|k_m}Q_{k_m|j_{m'}}\left[1+O(N^{-1/2})\right]
 \label{animalformula1}\ee
 if $h$ is a tree and 
\be O(N^{M-L})\label{animalformula2}\ .\ee
if $h$ has cycles. For trees, the $\ell$th edge has $m=\ell$, and $m'\in\{0\}\cup[\ell-1]$ numbers the node to which $w_\ell$ attaches.\end{theorem}

\begin{remarks} \begin{enumerate}
  \item 
Formula \eqref{animalformula2} shows clearly what is meant by saying that configuration graphs are {\it locally tree-like}\index{locally tree-like} as $N\to\infty$. It means the number of occurrences of any fixed finite size graph $h$ with cycles embedded within a configuration graph of size $N$ remains bounded with high probability as $N\to\infty$. 
  \item Even more interesting is that \eqref{animalformula1} shows that large configuration graphs exhibit a strict type of conditional independence. Selection of any root node $w_0$ of the tree graph $h$ splits it into two (possibly empty) trees $h_1, h_2$ with node-types $( j_{m}, k_{m}), m\in [M_1]$ and $( j_{m}, k_{m}),m \in [M_1+M_2]\setminus[M_1]$ where $M=M_1+M_2$. When we condition on the node-type of $w_0$, \eqref{animalformula1} shows that the remaining node-types form independent families: 
\beq \BBP[w_m\in\CN_{j_mk_m}, m\in[M], h\big | X, w_0\in\CN_{jk} ]&&\nonumber \\ &&\hspace{-2in}=\ \BBP[w_m\in\CN_{j_mk_m}, m\in[M_1], h_1\big |X, w_0\in\CN_{jk} ]\nonumber \\
 &&\hspace{-2in}\times \BBP[w_m\in\CN_{j_mk_m}, m\in[M_1+M_2]\setminus[M_1], h_2\big | X, w_0\in\CN_{jk} ] \ .\label{LTI_1}
\eeq
We call this deep property of the general configuration graph the {\it locally tree-like independence property}\index{locally tree-like independence property} (LTI property). In \cite{Hurd2015}, the LTI property provides the key to unravelling cascade dynamics in large configuration graphs. 
\end{enumerate} \end{remarks}

\bigskip\noindent{\bf Proof of Theorem \ref{animal_theorem}: \ } First, suppose Step 1 generates the node-type sequence $X$. Conditioned on $X$, now suppose the first step generates an in-edge $(w_1,w_0)$. Then, by refining Part 2 of Proposition \ref{P1},  the conditional probability that node $w_1$ has degree $j_1,k_1$ can be written
\beqq \frac{\BBP[w_1\in\CN_{j_1k_1},w_0\in\CN_{jk}|h,X]}{\BBP[w_0\in\CN_{jk}|h,X]}&&\\&&\hspace{-1.5in}=\ \frac{C^{-1}(e^-(0),e^+(0))e_{k_1,j_1}(0)e^-_{j,k}(0) Q_{k_1j}
C(e^-(1),e^+(1))}{C^{-1}(e^-(0),e^+(0))\sum_{k'}e^+_{k'}(0)e^-_{j,k}(0) Q_{k'j}C(e^-(1),e^+(1))}\\&&\hspace{-1.5in}=\ \left(\frac{e_{k_1,j_1}(0)e^-_{j,k}(0)}{e^+_{k_1}(0)e^-_{j}(0)}\right)\left(\frac{\BBE[e_{k_1j}|e^-(0), e^+(0)]}{E}\right)\left(\frac{e^-_{j,k}(0)}{E}\right)^{-1}\\
&&\hspace{-1.5in}=\ \left(\frac{k_1u_{k_1,j_1}}{k_1u^+_{k_1}}\right)\left(\frac{\BBE[e_{k_1j}|e^-(0), e^+(0)]}{e^-_{j}(0)}\right)\ .
\eeqq
Be aware that $C(e^-(1),e^+(1))$ in the denominator after the first equality depends on $k'$ and hence does not cancel a factor in the numerator.
Now, for $N\to\infty$, Part 2 of the Corollary applies to the second factor, and \eqref{LLN_P} applies to the first factor, and shows that for the case of an in-edge on the first step, with high probability, $X$ is such that:
\[  \BBP[w_1\in\CN_{j_1k_1}|w_0\in\CN_{jk},h,X] \ = \   P_{j_1|k_1}\ Q_{k_1|j}\left[1+O(N^{-1/2})\right]\ .
\]
The case of an out-edge is similar. 

Now we continue conditionally on $X$ from Step 1 and assume inductively that \eqref{animalformula1} is true for $M-1$ and prove it for $M$.  Suppose the final node $w_M$ is in-connected to the node $w_{M'}$ for some $M'\le M$. The ratio $\BBP[w_m\in\CN_{j_mk_m}, m\in[M]|v\in\CN_{jk},h,X]/\BBP[w_m\in\CN_{j_mk_m}, m\in[M-1]|w_0\in\CN_{jk},h,X]$ can be treated just as in the previous step and shown to be
\[\left(\frac{e_{k_M,j_M}(M-1)}{e^+_{k_M}(M-1)}\right)\left(\frac{\BBE[e_{k_Mj_{M'}}|e^-(M-1), e^+(M-1)]}{ e^-_{j_{M'}}(M-1)}\right)
\]
which with high probability equals \[  \BBP[w_1\in\CN_{j_1k_1}|w_0\in\CN_{jk},h,X] \ = \   P_{j_M|k_M}\ Q_{k_M|j_{M'}}\left[1+O(N^{-1/2})\right]\ .
\]
The case $w_M$ is out-connected to the node $w_{M'}$ is similar.

The first step $m$ that a cycle is formed can be treated by imposing a condition that $w_m=w_{m''}$ for some fixed $m''<m$. One finds that the conditional probability of this is 
\beqq \BBP[w_m=w_{m''}, w_\ell \in\CN_{j_\ell k_\ell}, \ell\in[m-1]|w_0\in\CN_{jk},h,X]&&\\
&&\hspace{-2.5in}=\frac{k_{m''}}{{e^+_{k_{m''}}(m-1)}}\ \times\ \BBP[w_\ell \in\CN_{j_\ell k_\ell}, \ell\in[m-1]|w_0\in\CN_{jk},h,X]\ .
\eeqq
The first factor is $O(N^{-1})$ as $N\to\infty$, which proves the desired statement \eqref{animalformula2} for cycles. 

Finally, since \eqref{animalformula2} is true for cycles, with high probability all finite configurations are trees. Therefore their asymptotic probability laws are given by \eqref{animalformula1}, as required.

\qedsymbol

\bigskip

\section{Approximate ACG Simulation} \label{sec:5}

It was observed in Section \ref{Configuration} that Step 1 of the configuration graph construction draws a sequence $(j_i,k_i)_{i\in[N]}$ of node types that is iid with the correct distribution $P$, but is only feasible, $\sum_i(k_i-j_i)=0$, with small probability. Step 2 of the exact ACG algorithm in Section \ref{finite_assortative_section} requires is even less feasible in practice. Practical simulation algorithms address the first problem by ``clipping'' the drawn node bidegree sequence when the discrepancy $D=D_N:=\sum_i(k_i-j_i)$ is not too large, meaning it is adjusted by a small amount to make it feasible, without making a large change in the joint distribution. Step 1 of the following simulation algorithm generalizes slightly the method introduced by \cite{ChenOlve13} who verify that the effect of clipping vanishes with high probability as $N\to\infty$. The difficulty with Step 2 of the ACG construction is overcome by an approximate sequential wiring algorithm. 

The {\it approximate assortative configuration simulation algorithm} for multigraphs of size $N$, parametrized by the node-edge degree distribution pair $(P,Q)$ that have support on the finite set $(j,k)\in\{0,1,\dots, K\}^2$,  involves choosing a suitable threshold $T=T(N)$ and modifying the steps identified in Section \ref{finite_assortative_section}: 
\begin{enumerate}
  \item Draw a sequence of  $N$ node-type pairs $X= ((j_1,k_1),\dots, (j_N,k_N))$ independently from $P$, and accept the draw if and only if $0<|D|\le T(N)$. When the sequence $(j_i,k_i)_{i\in[N]}$ is accepted,   the sequence is adjusted by adding a few stubs, either in- or out- as needed. First draw a random subset $\sigma\subset\CN$ of size $|D|$ with uniform probability $\binom{N}{|D|}^{-1}$, and then define the feasible sequence $\tilde X=(\tilde j_i,\tilde k_i)_{i\in[N]}$ by adjusting the degree types for $i\in\sigma$ as follows:
\beq
\label{clipone}
\tilde j_i= j_i + \xi^-_i;\quad \xi^-_i&=&\One(i\in\sigma,D>0)\\
\label{cliptwo}
\tilde k_i= k_i + \xi^+ _i;\quad \xi^+_i&=&\One(i\in\sigma,D<0)\ .
\eeq
   \item Conditioned on $\tilde X$, the result of Step 1, randomly wire together available in and out stubs {\it sequentially}, with suitable weights, to produce the sequence of edges $W$. At each $\ell=1,2,\dots, E$, match from available in-stubs and  out-stubs weighted according to their degrees $j,k$ by 
\be\label{ACWweights}
C^{-1}(\ell)\frac {Q_{kj}}{Q^+_kQ^-_j}\ .
\ee In terms of the bivariate random process $( e^-_j(\ell),  e^+_k(\ell))$ with initial values\\ $( e^-_j(1), e^+_k(1))=(e^-_j, e^+_k)$  that at each $\ell$ counts the number of available degree $j$ in-stubs and  degree $k$ out-stubs, the $\ell$ dependent normalization factor $C(\ell)$ is given by: \be
\label{Cfactor}
C(\ell)=\sum_{jk}  e^-_j(\ell)  e^+_k(\ell) \frac {Q_{kj}}{Q^+_kQ^-_j}
\ .
\ee  
\end{enumerate}

\begin{remark} An alternative simulation algorithm for the ACG model has been proposed and studied in \cite{DeprWuth15}.
\end{remark}

Chen and Olvera-Cravioto, \cite{ChenOlve13}, addresses the clipping in Step 1 and shows that the discrepancy of the approximation is negligible as $N\to\infty$:
\begin{theorem}
\label{nodedegreetheorem} Fix $\delta\in(0,1/2)$, and for each $N$ let the threshold be $T(N)=N^{1/2+\delta}$. Then:\begin{enumerate}
 \item The acceptance probability $ \BBP[|D_N|\le T(N)]\to 1$
as $N\to\infty$;
 \item 
 For any fixed finite $M$, $\Lambda$, and bounded function $f:(\BBZ_+\times{\BBZ_+})^{M}\to[-\Lambda,\Lambda]$
\be \left| \BBE[f\left((\tilde j_i,\tilde k_i)_{i=1,\dots, M}\right)]-\BBE[f\left((\hat j_i,  \hat k_i)_{i=1,\dots,M}\right)]\right|\to 0\ ;
\ee
\item The following limits in probability hold:
\be \frac1{N}\tilde u_{jk}\Parrow P_{jk},\qquad \frac1{N}\tilde u^+_k  \Parrow P^+_{k},\qquad \frac1{N}\tilde u^-_{j}\Parrow P^-_{j}\ .\label{empirical}
\ee
 \end{enumerate}
 \end{theorem}

Similarly it is intuitively clear that the discrepancy of the approximation in Step 2 is negligible as $N\to\infty$. As long as $e^-_j(\ell),  e^+_k(\ell)$ are good approximations of $(E-\ell)Q^-_j,  (E-\ell)Q^+_k$, \eqref{ACWweights} shows that the probability that edge $\ell$ has type $(k,j)$ will be approximately $Q_{kj}$. Since the detailed analysis of this problem is not yet complete, we state the desired properties as a conjecture:
\begin{conjecture}
\label{Approx_assortative} In the approximate assortative configuration graph construction with probabilities $P,Q$, the following convergence properties hold as $N\to\infty$. 
\begin{enumerate}
  \item The fraction  of type $(k,j)$ edges in the matching sequence $( k_{l}, j_{\ell})_{\ell\in[E]}$ concentrates with high probability around the nominal edge distribution $Q_{kj}$:
\be 
\frac{{e}_{kj}}{E} = Q_{kj}+o(1)\ .
\ee
  \item For any fixed finite number $L$, the first $L$ edges $\ell,\ell\in[L]$ have degree sequence $(k_{l},j_{\ell})_{\ell\in[L]}$ that converges in distribution to $(\hat k_{l},\hat j_{\ell})_{\ell\in[L]}$, an independent  sequence of identical $Q$ distributed random variables.
\end{enumerate}\end{conjecture}

Although the conjecture is not yet completely proven, extensive simulations have verified the consistency of the approximate ACG algorithm with the theoretical large $N$ probabilities.

\bibliographystyle{plain}


\begin{thebibliography}{10}

\bibitem{BechAtal10}
M.~Bech and E.~Atalay.
\newblock The topology of the federal funds market.
\newblock {\em Physica A: Statistical Mechanics and its Applications},
  389(22):5223--5246, 2010.

\bibitem{Bollobas80}
B.~Bollob\~as.
\newblock A probabilistic proof of an asymptotic formula for the number of
  labelled regular graphs.
\newblock {\em Eur. J. Comb.}, 1:311, 1980.

\bibitem{Bollobas01}
B.~Bollob\~as.
\newblock {\em Random Graphs}.
\newblock Cambridge studies in advanced mathematics. Cambridge University
  Press, 2 edition, 2001.

\bibitem{ChenOlve13}
Ningyuan Chen and Mariana Olvera-Cravioto.
\newblock Directed random graphs with given degree distributions.
\newblock {\em Stochastic Systems}, 3(1):147--186, 2013.

\bibitem{DeprWuth15}
Philippe Deprez and Mario~V. W\"uthrich.
\newblock Construction of directed assortative configuration graphs.
\newblock arXiv:1510.00575, October 2015.

\bibitem{Erdelyi1956}
Arthur Erd{\'e}lyi.
\newblock {\em Asymptotic expansions}.
\newblock Dover, New York, 1956.

\bibitem{ErdoReny59}
P.~Erd\"os and A.~R\'enyi.
\newblock On random graphs.
\newblock {\em I. Publ. Math. Debrecen}, 6:290--297, 1959.

\bibitem{Hurd10}
T.~R. Hurd.
\newblock Saddlepoint approximation.
\newblock In Rama Cont, editor, {\em Encyclopedia of Quantitative Finance}.
  John Wiley \& Sons, Ltd, 2010.

\bibitem{Hurd2015}
T.~R. Hurd.
\newblock Contagion! systemic risk in financial networks.
\newblock Available at http://ms.mcmaster.ca/tom/tom.html, 2015.

\bibitem{Janson09_simple}
Svante Janson.
\newblock The probability that a random multigraph is simple.
\newblock {\em Combinatorics, Probability and Computing}, 18:205--225, 3 2009.

\bibitem{MayArin10}
Robert~M. May and Nimalan Arinaminpathy.
\newblock Systemic risk: the dynamics of model banking systems.
\newblock {\em Journal of The Royal Society Interface}, 7(46):823--838, 2010.

\bibitem{MollReed95}
Michael Molloy and Bruce Reed.
\newblock A critical point for random graphs with a given degree sequence.
\newblock {\em Random Structures \& Algorithms}, 6(2-3):161--180, 1995.

\bibitem{Soramakietal07}
Kimmo Soram{\"a}ki, M.~Bech, J.~Arnold, R.~Glass, and W.~Beyeler.
\newblock The topology of interbank payment flows.
\newblock {\em Physica A: Statistical Mechanics and its Applications},
  379(1):317--333, 2007.

\bibitem{vdHofstad14}
R.~van~der Hofstad.
\newblock Random graphs and complex networks.
\newblock Book, to be published, 2014.

\end{thebibliography}

\end{document}